\documentclass[twoside]{article}
\usepackage{graphicx}
\usepackage{amsmath, amsthm, amssymb}

\newtheorem{theorem}{Theorem}
\newtheorem{remark}{Remark}
\newtheorem{lemma}{Lemma}
\newtheorem{corollary}{Corollary}
\def \hs {\hspace*{-2mm}}
\newfont{\mi}{cmti9}
\newfont{\m}{cmr8}
\newfont{\ms}{cmsl8}
\newfont{\autor}{cmcsc10}
\begin{document}
\rule[3mm]{128mm}{0mm}
\vspace*{-16mm}

{\footnotesize K\,Y\,B\,E\,R\,N\,E\,T\,I\,K\,A\, ---
\,V\,O\,L\,U\,M\,E\, {\it 4\,4\,} (\,2\,0\,0\,8\,)\,,\,
N\,U\,M\,B\,E\,R\, X\,,\,\  P\,A\,G\,E\,S\, \,x\,x\,x -- x\,x\,x}\\
\rule[3mm]{128mm}{0.2mm}

\vspace*{11mm}

{\large\bf \noindent Optimal  sequential multiple hypothesis tests}

\vspace*{8mm}

{\autor \indent Andrey Novikov}\vspace*{23mm}\\
\small This work deals with a general problem of testing multiple
hypotheses about the distribution of a discrete-time stochastic
process. Both the Bayesian and the conditional settings are
considered. The structure
 of optimal sequential tests is characterized.

\smallskip\par
\noindent {\sl Keywords:}\,
\begin{minipage}[t]{112mm}
 sequential analysis,
hypothesis testing, multiple hypotheses, discrete-time stochastic
process, dependent observations, optimal
 sequential test, Bayes sequential test
\end{minipage} \smallskip

\normalsize

\section{\hs INTRODUCTION}\label{s1}
Let $X_1,X_2,\dots, X_n, \dots$ be a discrete-time stochastic
process, whose distribution depends on   an unknown "parameter"
$\theta$. We consider the classical problem of testing multiple
hypotheses $H_1:\, \theta=\theta_1$, $H_2:\, \theta=\theta_2$,
$\dots$, $H_k:\, \theta=\theta_k$, $k\geq 2$.

The main goal of this article is to characterize the structure of
optimal sequential tests in this problem.

Let us  suppose that for any $n=1,2,\dots,$ the vector $(X_1, X_2,
\dots , X_n)$ has a probability "density" function
\begin{equation}\label{0}
 f_\theta^{n}(x_1, x_2, \dots, x_n)
 \end{equation}
 (Radon-Nicodym derivative of its distribution) with respect to
 a product-measure $$\begin{array}{cc}
                            \mu^n=& \underbrace{\mu\otimes \mu\otimes \dots\otimes\mu}, \\
                             &n \;\mbox{times}
                           \end{array}
 $$
 for some $\sigma$-finite measure $\mu$ on the respective space.

 We define a (randomized) sequential hypothesis test
 as a pair $(\psi, \phi)$ of a {\em stopping
rule} $\psi$ and a {\em decision rule} $\phi$,
 with $$\psi=\left(\psi_1,\psi_2, \dots ,\psi_{n},\dots\right),$$
and
$$\phi=\left(\phi_1,\phi_2, \dots ,\phi_n,\dots\right).$$
The  functions $$\psi_n=\psi_n(x_1, x_2,\dots, x_n),\quad
n=1,2,\dots,$$ are supposed to be some measurable functions with
values in $[0, 1]$. The functions
$$\phi_n=\phi_n(x_1, x_2,\dots, x_n),\quad n=1,2,\dots$$
are supposed to be measurable vector-functions with $k$ non-negative
components $\phi_n^i=\phi_n^i(x_1,\dots,x_n)$:
$$
\phi_n=(\phi_n^1,\dots,\phi_n^k),
$$
such that $\sum_{i=1}^k \phi_n^i=1$ for any $n=1,2,\dots$.

The interpretation of all these elements is as follows.

The value of  $\psi_n(x_1,\dots,x_n)$ is interpreted  as the
conditional probability {\em to stop and proceed to decision
making}, given that we came to stage $n$ of the experiment and that
the observations up to stage $n$ were $(x_1, x_2, \dots, x_n).$ If
there is no stop, the experiments continues to the next stage and an
additional observation $x_{n+1}$ is taken. Then the rule
$\psi_{n+1}$ is applied to $x_1,x_2,\dots, x_n,x_{n+1}$ in the same
way as as above, etc., until the experiment eventually stops.

It is supposed that when the experiment stops,  a decision {\em to
accept some of $H_1,\dots, H_k$} is to be made. The function
$\phi_n^i(x_1,\dots, x_n)$ is interpreted  as the conditional
probability {\em to accept} $H_i$, $i=1,\dots,k$, given that the
experiment stops at stage $n$ being  $(x_1,\dots, x_n)$ the data
vector observed.

The stopping rule $\psi$ generates, by the above process, a random
variable $\tau_\psi$ ({\em stopping time}) whose distribution is
given by
$$
P_\theta(\tau_\psi=n)=E_\theta (1-\psi_1)(1-\psi_2)\dots
(1-\psi_{n-1})\psi_n.
$$
Here, and throughout the paper, we interchangeably use  $\psi_n$
both for $\psi_n(x_1,x_1,\dots,x_n)$ and for
$\psi_n(X_1,X_1,\dots,X_n)$, and so do we for any other function of
observations $F_n$. This does not cause any problem if we adopt the
following agreement: when $F_n$ is under probability or expectation
sign, it is $F_n(X_1,\dots, X_n)$, otherwise it is $F_n(x_1,\dots,
x_n)$.

For a sequential test $(\psi,\phi)$ let us define
\begin{equation}\label{3a}
\alpha_{ij}(\psi,\phi)=P_{\theta_i}(\,\mbox{accept}\,
H_j)=\sum_{n=1}^\infty E_{\theta_i}(1-\psi_1)
\dots(1-\psi_{n-1})\psi_n\phi_{n}^j
\end{equation}
and
\begin{equation}\label{3c}
\beta_{i}(\psi,\phi)=P_{\theta_i}(\,\mbox{accept any}\,
H_j\,\mbox{different from}\,
H_i)=\sum_{j\not=i}\alpha_{ij}(\psi,\phi),
\end{equation}
$i=1,\dots,k$, $j=1,\dots, k$. The probabilities
$\alpha_{ij}(\psi,\phi)$ for $j\not=i$ can be considered
"individual" error probabilities and $\beta_i(\psi,\phi)$ "gross"
error probability, under hypothesis $H_i$, of the sequential test
$(\psi,\phi)$.

Another important characteristic of a sequential test is the {\em
average sample number}:
\begin{equation}\label{8aa}
N(\theta;\psi)=E_\theta\tau_\psi=\begin{cases}\sum_{n=1}^\infty
nP_\theta(\tau_\psi=n), \;\mbox{if}\;
P_\theta(\tau_\psi<\infty)=1,\cr
\infty\quad\mbox{otherwise}.\end{cases}
\end{equation}

Let $\theta$ be any fixed (and known)
value of the parameter (we {\em do not}
suppose, generally, that $\theta$ is
one of $\theta_i$, $i=1,\dots k$).

In this article, we solve the two
following problems:
\begin{description}
\item[Problem I.]
Minimize $N(\psi)=N(\theta;\psi)$  over
all sequential tests $(\psi,\phi)$
subject to
\begin{equation}\label{1}
\alpha_{ij}(\psi,\phi)\leq\alpha_{ij},\quad\mbox{for all}\;
i=1,\dots k,\;\mbox{and for all}\; j\not= i,
\end{equation}
where $\alpha_{ij}\in(0,1)$ (with $i,j=1,\dots k$, $j\not=i$) are
some constants.
\item[Problem II.]
Minimize $N(\psi)$  over all sequential
tests $(\psi,\phi)$ subject to
\begin{equation}\label{2}
\beta_{i}(\psi,\phi)\leq\beta_i,\quad\mbox{for all}\; i=1,\dots
k,
\end{equation}
with some constants $\beta_i\in(0,1)$, $i=1,\dots,k$.
\end{description}

More general problems of minimizing an average cost of type
$$
N(\psi)=\sum_{n=1}^\infty E_\theta C_n
(1-\psi_1) \dots(1-\psi_{n-1})\psi_n
$$
with some cost function $C_n=C_n(X_1, X_2,\dots X_n)$ can be treated
in essentially the same manner.

If $k=2$ then Problems I and II are equivalent, because
$\beta_1(\psi,\phi)=\alpha_{12}(\psi,\phi)$, and
$\beta_2(\psi,\phi)=\alpha_{21}(\psi,\phi)$, by (\ref{3a}) and
(\ref{3c}).

For independent and identically
distributed (i.i.d.) observations and
$k=2$ the formulated problem, when
$\theta\not =\theta_1$ and $\theta\not
=\theta_2$, is known as the modified
Kiefer-Weiss problem (see
\cite{Weiss}), being the original
Kiefer-Weiss problem minimizing
$\sup_\theta N(\psi)$ under (\ref{1})
(see \cite{KieferWeiss}).

For the latter problem,  taking into
account the usual relations between
Bayesian and minimax procedures, it
seems to be reasonable to generalize
our problem of minimizing
$N(\theta;\psi)$ to that of minimizing
$$
\int N(\theta;\psi)d\pi(\theta),
$$
with some "weight"
measure $\pi$.  From what follows it is easily
seen that, under natural measurability  conditions,
our method works as well for this latter problem.

In Section 2, we reduce Problems I and II  to an unconstrained
minimization problem. The new objective function is the
Lagrange-multiplier function $L(\psi;\phi)$.

In Section 3, we find
$$L(\psi)=\inf_{\phi}L(\psi,\phi),$$ where the infimum is taken over
all decision rules.

In Section 4, we minimize $L(\psi)$ in the class of truncated
stopping rules, i.e. such that $\psi_N\equiv 1$ for some
$0<N<\infty$.

In Section 5, we characterize the structure of optimal stopping rule
$\psi$ in the class of all stopping rules.

In Section 6, we apply the results obtained in Sections 2 -- 5 to
the solution of Problems I and II.

\section{\hs REDUCTION TO NON-CONSTRAINED
MINIMIZATION}\label{s2}

In this section, the Problems I and II will be reduced to
unconstrained optimization problems using the idea of the Lagrange
multipliers method.

\subsection{\normalsize Reduction to Non-Constrained Minimization in Problem I}
 To proceed with minimizing $N(\psi)$ over the sequential tests
subject to (\ref{1}), let us define the following
Lagrange-multiplier function:
\begin{equation}\label{4}
L(\psi,\phi)=N(\psi)+\sum_{{1\leq
i,j\leq
k};\,{i\not=j}}\lambda_{ij}\alpha_{ij}(\psi,\phi)
\end{equation}
where  $\lambda_{ij}\geq 0$   are some constant multipliers. Recall
that $N(\psi)=E_\theta\tau_\psi$, where $\theta$ is the fixed value
of parameter for which the average sample number (\ref{8aa}) is to
be minimized. Generally, we {\em do not} suppose that $\theta$ is
one of $\theta_i$, $i=1,\dots, k$.

Let $\Delta$ be a class of tests.

The following theorem is a direct application of the Lagrange
multipliers method.

\begin{theorem}\label{t1} Let exist $\lambda_{ij}> 0$, $i=1,\dots, k$, $j=1,\dots, k$, $j\not=i$, and a test $(\psi^*,\phi^*)\in
\Delta$  such that for all sequential tests  $(\psi,\phi)\in
\Delta$
\begin{equation}\label{5}
L(\psi^*,\phi^*)\leq L(\psi,\phi)
 \end{equation}
 holds and such that
 \begin{equation}\label{6}\alpha_{ij}(\psi^*,\phi^*)=\alpha_{ij}\quad\mbox{for
 all}\quad
i=1,\dots k,\;\mbox{and for all}\; j\not= i.
 \end{equation}

 Then for all $(\psi,\phi)\in\Delta$ such
 that
 \begin{equation}\label{5bis}
 \alpha_{ij}(\psi,\phi)\leq\alpha_{ij}\quad\mbox{for all}\quad
i=1,\dots k,\;\mbox{and for all}\; j\not= i,
 \end{equation}
 it holds
\begin{equation}\label{5a}
N(\psi^*)\leq  N(\psi).
\end{equation}

 The inequality
in (\ref{5a}) is strict if at least one of the equalities
(\ref{5bis}) is strict.
\end{theorem}

\begin{proof}  Let $(\psi,\phi)\in \Delta$ be any sequential test
satisfying (\ref{5bis}). Because of (\ref{5})
\begin{eqnarray}\nonumber L(\psi^*,\phi^*)&=&N(\psi^*)+\sum_{j\not=
i}\lambda_{ij}\alpha_{ij}(\psi^*,\phi^*)\\
\label{5c}&\leq&
L(\psi,\phi)=N(\psi)+\sum_{j\not=i}\lambda_{ij}\alpha_{ij}(\psi,\phi)\leq N(\psi)+\sum_{j\not =
i}\lambda_{ij}\alpha_{ij},
\end{eqnarray}
where to get the last inequality we used (\ref{1}).

So,
$$
N(\psi^*)+\sum_{j\not=
i}\lambda_{ij}\alpha_{ij}(\psi^*,\phi^*)
\leq
N(\psi)+\sum_{j\not=i}\lambda_{ij}\alpha_{ij},
$$
and taking into account conditions (\ref{6}) we get from this that
$$N(\psi^*)\leq N(\psi).
$$

To get the last statement of the
theorem we note that if
$N(\psi^*)=N(\psi)$ then there are
equalities in (\ref{5c})
instead of inequalities which is only
possible if
$\alpha_{ij}(\psi,\phi)=\alpha_{ij}$
for any $i,j=1,\dots k$, $j\not=i$.
\end{proof}
\begin{remark}\label{r0}
\rm The author owes the idea of the use of the Lagrange-multiplier
method in sequential hypotheses testing to Berk \cite{Berk}.
Essentially, the method of Lagrange multipliers is implicitly used
in the monograph of Lehmann \cite{Lehmann} in the proof of  the
fundamental lemma of Neyman-Pearson. In a way, the Bayesian approach
in hypotheses testing can be considered as a variant of the
Lagrange-multiplier method as well.
\end{remark}
\begin{remark}\label{r00}\rm
All our results below can be adapted to the Bayesian context by
choosing appropriate Lagrange multipliers and using
$$\sum_{i=1}^k N(\theta_i;\psi)\pi_i$$
instead of $N(\theta;\psi)$ in
$L(\psi,\phi)$ above. From this point
of view, we extend and complement the
results of Cochlar \cite{Cochlar} about
the { existence} of Bayesian
sequential tests.

More generally, all our results are applicable as well for minimization of
$$\int N(\theta;\psi)d\pi(\theta),$$ where $\pi$ is any probability measure (see Remarks \ref{r5} and \ref{r7a} below).
\end{remark}
\subsection{\normalsize Reduction to Non-Constrained Minimization in Problem II}

Very much like in the preceding section, define
\begin{equation}\label{4aBis}
L(\psi,\phi)=N(\psi)+\sum_{i=1}^k\lambda_{i}\beta_{i}(\psi,\phi),
\end{equation}
where $\lambda_i\geq 0$ are the Lagrange multipliers.

In a very similar manner to Theorem \ref{t1}, we have
\begin{theorem}\label{t1a} Let exist $\lambda_{i}> 0$, $i=1,\dots, k$,  and a sequential test $(\psi^*,\phi^*)\in
\Delta$ such that for all $(\psi,\phi)\in \Delta$
\begin{equation}\label{5Bis}
L(\psi^*,\phi^*)\leq L(\psi,\phi)
 \end{equation}
 holds and such that
 \begin{equation}\label{6BisBis}\beta_{i}(\psi^*,\phi^*)=\beta_{i}\quad\mbox{for
 all}\quad
i=1,\dots k.
 \end{equation}

 Then for all sequential tests $(\psi,\phi)\in\Delta$ such
 that
 \begin{equation}\label{5bisBis}
 \beta_{i}(\psi,\phi)\leq\beta_{i}\quad\mbox{for all}\quad
i=1,\dots k,
 \end{equation}
 it holds
\begin{equation}\label{5aBis}
N(\psi^*)\leq  N(\psi).
\end{equation}

 The inequality
in (\ref{5aBis}) is strict if at least one of the equalities
(\ref{5bisBis}) is strict.
\end{theorem}

\section{\hs OPTIMAL DECISION RULES}\label{s3}

Due to Theorems \ref{t1} and \ref{t1a}, Problem I is reduced to
minimizing (\ref{4}) and Problem II is reduced to minimizing
(\ref{4aBis}). But (\ref{4aBis}) is a particular case of (\ref{4}),
namely, when $\lambda_{ij}=\lambda_i$ for any $j=1,\dots,k$,
$j\not=i$ (see (\ref{3a}) and (\ref{3c})). Because of that, from now
on, we will only solve the problem of minimizing $L(\psi,\phi)$
defined by (\ref{4}).

In particular, in this section  we  find
$$
L(\psi)=\inf_{\phi}L(\psi,\phi),
$$
and the corresponding decision rule $\phi$, at which this infimum is
attained.

 Let $I_A$ be the
indicator function of the event $A$.

\begin{theorem} \label{t2} For any $\lambda_{ij}\geq 0$, $i=1,\dots,k$, $j\not= i$,  and for
any sequential test $(\psi,\phi)$
\begin{equation}\label{7}
L(\psi,\phi)\geq
N(\psi)+\sum_{n=1}^\infty \int
(1-\psi_1)\dots (1-\psi_{n-1})\psi_n
l_n d\mu^n,
\end{equation}
where
\begin{equation}\label{6a}
l_n=\min_{1\leq j\leq k}\sum_{i\not =j}\lambda_{ij}f_{\theta_i}^n.
\end{equation}
Supposing that $N(\psi)$ is finite, the
right-hand side of (\ref{7}) is
attained if and only if
\begin{equation}\label{7aa}
\phi_{n}^j\leq I_{\left\{\sum_{i\not
=j}\lambda_{ij}f_{\theta_i}^n=l_n\right\}}
\end{equation}
for all $j=1,\dots k$, $\mu^n$-almost anywhere on
$$S_n^\psi=\{(x_1,\dots,x_n): s_n(x_1,\dots,x_n)>0\},$$
where $s_n^\psi(x_1,\dots,x_n)=s_n^\psi=(1-\psi_1)\dots(1-\psi_{n-1})\psi_n$, for all $n=1,2,\dots$.
\end{theorem}
\begin{proof} Inequality (\ref{7}) is equivalent to
\begin{equation}\label{9}
\sum_{1\leq i,j\leq
k;\,j\not=i}\lambda_{ij}\alpha_{ij}(\psi,\phi)\geq\sum_{n=1}^\infty
\int (1-\psi_1)\dots (1-\psi_{n-1})\psi_n l_n d\mu^n.
\end{equation}

We prove it by finding a lower bound for the left-hand side of
(\ref{9}) and proving that this lower bound is attained
 if $\phi$ satisfies  (\ref{7aa}).

To do this, we will use the following simple
\begin{lemma}\label{l0}
Let $\phi_1,\dots,\phi_k$ and  $F_1, \dots F_k$ be some measurable
non-negative functions on a measurable space with a measure $\mu$,
such that
$$
 \sum_{i=1}^k \phi_i(x)\equiv 1,
$$
and such that
$$
\int \min_{1\leq i\leq k}F_i(x)d\mu(x)<\infty.
$$
 Then
\begin{equation}\label{107b}
\int \left(\sum_{i=1}^k\phi_i(x)F_i(x)\right)d\mu(x)\geq \int
\min_{1\leq i\leq k}F_i(x)d\mu(x)
\end{equation}
with an equality in (\ref{107b}) if and only if
\begin{equation}\label{107e}
\phi_i\leq I_{\displaystyle{\{F_i= \min_{1\leq j\leq k}
F_j\}}}\;\mbox{ for any}\quad i=1,2,\dots,k,
\end{equation}
$\mu$-almost anywhere.

\end{lemma}
\begin{proof} To prove (\ref{107b}) it suffices to show that
\begin{equation}\label{107c}
\int \left(\sum_{i=1}^k\phi_i(x)F_i(x)\right)d\mu(x)-\int \min_{1\leq i\leq k}F_i(x)d\mu(x)\geq 0,
\end{equation}
because the second integral is finite by the conditions of the
Lemma.

But (\ref{107c}) is equivalent to
\begin{equation}\label{107d}
\int \sum_{i=1}^k\phi_i(x)(F_i(x)-\min_{1\leq j\leq
k}F_j(x))d\mu(x)\geq 0,
\end{equation}
being this trivial because the function under the integral sign is
non-negative.

Because of  this, there is an equality in (\ref{107d}) if and only
if
$$
\sum_{i=1}^k\phi_i(x)(F_i(x)-\min_{1\leq j\leq k} F_j(x))=0
$$
$\mu$-almost anywhere, which is only possible if  (\ref{107e}) holds
true. \end{proof}

Starting with the proof of (\ref{9}), let us give to the left-hand
side of it
 the form
 $$
\sum_{1\leq i,j\leq k;\,j\not=i}\lambda_{ij}\alpha_{ij}(\psi,\phi)
 $$
\begin{equation}\label{7f}
= \sum_{n=1}^\infty \int (1-\psi_1)\dots (1-\psi_{n-1})\psi_n
 \sum_{j=1}^k\left(\sum_{1\leq i\leq
k;\,i\not=j}\lambda_{ij}f_{\theta_i}^{n}\right)\phi_{n}^j d\mu^n
 \end{equation}
(see (\ref{3a})).

Applying Lemma \ref{l0} to each summand in (\ref{7f}) we immediately
have:
 \begin{equation}\label{7g}
\sum_{1\leq i,j\leq k;\,j\not=i}\lambda_{ij}\alpha_{ij}(\psi,\phi)
\geq \sum_{n=1}^\infty \int (1-\psi_1)\dots (1-\psi_{n-1})\psi_n
 l_n
d\mu^n
 \end{equation}
with an equality if and only if
$$
\phi_{n}^j\leq
I_{\{\sum_{i\not=j}\lambda_{ij}f_{\theta_i}^{n}=l_n\}}
$$
for all $1\leq j\leq k$, $\mu^n$-almost anywhere on $S_n^\psi$,
 for all
$n=1,2,\dots$.\end{proof}
\begin{remark}\label{r2}  It is easy to see, using (\ref{8aa}) and (\ref{7g}), that
\begin{equation}\label{10}
L(\psi)=\inf_{\phi}L(\psi,\phi)=\sum_{n=1}^\infty
\int (1-\psi_1)\dots
(1-\psi_{n-1})\psi_n\left(nf_{\theta}^{n}+l_n
\right) d\mu^n
\end{equation}
if $P_{\theta}(\tau_\psi<\infty)=1$ and $L(\psi)=\infty$ otherwise.
\end{remark}
\begin{remark}\label{r8}
In the Bayesian context of Remark \ref{r00}, the "if"-part of
Theorem \ref{t2} can also be  derived  from Theorem 5.2.1 \cite{Ghosh}.
\end{remark}

\section{\hs TRUNCATED STOPPING RULES}\label{s4}
Our next goal is to find a stopping rule $\psi$ minimizing the value
of $L(\psi)$ in (\ref{10}).

In this section, we solve, as an intermediate step,  the problem of
minimization of $L(\psi)$ in the class of truncated stopping rules,
that is, in the class $\Delta^N$ of
\begin{equation}\label{11}
\psi=(\psi_1,\psi_2,\dots,\psi_{N-1},1,\dots).\end{equation}

For any $\psi\in\Delta^N$ let us define
\begin{equation}\label{14}
L_N(\psi)=\sum_{n=1}^N \int (1-\psi_1)\dots
(1-\psi_{n-1})\psi_n\left(nf_{\theta}^{n}+l_n \right) d\mu^n
\end{equation} (see (\ref{10})).

 The following lemma takes over a large part of
work of minimizing $L_N(\psi)$ over $\psi\in\Delta^N$.

\begin{lemma}\label{l2} Let $r\geq 2$ be any natural number, and let $v_r=v_r(x_1,x_2,\dots,x_r)$ be
any measurable function such that $\int v_rd\mu^r<\infty$. Then
$$
\sum_{n=1}^{r-1}\int
(1-\psi_1)\dots(1-\psi_{n-1})\psi_n(nf_\theta^n+l_n)d\mu^n
$$
$$
+\int
(1-\psi_1)\dots(1-\psi_{r-1})\left(rf_\theta^r+v_r\right)d\mu^r
$$
\begin{equation}\label{37}
\geq\sum_{n=1}^{r-2}\int
(1-\psi_1)\dots(1-\psi_{n-1})\psi_n(nf_\theta^n+l_n)d\mu^n
\end{equation}
$$
+\int
(1-\psi_1)\dots(1-\psi_{r-2})\left((r-1)f_\theta^{r-1}+v_{r-1}\right)d\mu^{r-1},
$$
where
\begin{equation}\label{38}
v_{r-1}=\min\{l_{r-1},f_\theta^{r-1}+\int
v_r(x_1,\dots,x_r)d\mu(x_r)\}.
\end{equation}
There is an equality in (\ref{37}) if and only if
\begin{equation}\label{39}
I_{\{l_{r-1}< f_\theta^{r-1}+\int v_r(x_1,\dots,x_r)d\mu(x_r)\}}\leq
\psi_{r-1}\leq I_{\{l_{r-1}\leq f_\theta^{r-1}+\int
v_r(x_1,\dots,x_r)d\mu(x_r)\}},
\end{equation}
$\mu^{r-1}$-almost anywhere on $C_{r-1}^\psi$, where, by definition,
$$C_n^\psi=\{(x_1,\dots,x_n):(1-\psi_1(x_1))\dots(1-\psi_{n-1}(x_1,\dots,x_{n-1}))>0\}$$
for any $n=1,2,\dots$
\end{lemma}
\begin{proof}
Let us start with the following simple consequence of  Lemma
\ref{l0}.
\begin{lemma} \label{l1} Let $\chi,\phi, F_1, F_2$ be some measurable functions on
a measurable space with a measure $\mu$, such that
$$
0\leq\chi(x)\leq 1,\quad 0\leq\phi(x)\leq 1, \quad F_1(x)\geq 0,
\quad F_2(x)\geq 0,\quad
$$
and
$$
\int \min\{F_1(x),F_2(x)\}d\mu(x)<\infty.
$$

Then
\begin{equation}\label{7b}
\int \chi(x)(\phi(x)F_1(x)+(1-\phi(x))F_2(x))d\mu(x)\geq \int
\chi(x)\min\{F_1(x),F_2(x)\}d\mu(x)
\end{equation}
with an equality if and only if
\begin{equation}\label{7e}
I_{\{F_1(x)<F_2(x)\}}\leq \phi(x)\leq I_{\{F_1(x)\leq F_2(x)\}}
\end{equation}
$\mu$-almost anywhere on $\{x:\chi(x)>0\}$.
\end{lemma}
\begin{proof}
Defining $\phi_1(x)\equiv \phi(x)$ and $\phi_2(x)\equiv 1-\phi(x)$,
from Lemma \ref{l0} we immediately obtain (\ref{7b}), with an
equality if and only if
\begin{equation}\label{12}
\phi_1(x)=\phi(x)\leq
I_{\{\chi(x)(F_1(x)-\min\{F_1(x),F_2(x)\})=0\}}
\end{equation}
and
\begin{equation}\label{13}
\phi_2(x)=1-\phi(x)\leq
I_{\{\chi(x)(F_2(x)-\min\{F_1(x),F_2(x)\})=0\}}
\end{equation}
$\mu$-almost anywhere. Expressing $\phi(x)$ from (\ref{12}) and
(\ref{13}) we have that there is an equality in (\ref{7b}) if and
only if
\begin{equation*}
I_{\{\chi(x)(F_2(x)-\min\{F_1(x),F_2(x)\})>0\}}\leq \phi(x)\leq
I_{\{\chi(x)(F_1(x)-\min\{F_1(x),F_2(x)\})=0\}}
\end{equation*}
$\mu$-almost anywhere, which is equivalent to
\begin{equation*}
I_{\{F_1(x)<F_2(x)\}}\leq \phi(x)\leq I_{\{F_1(x)\leq F_2(x)\}}\quad
\mu-\mbox{almost anywhere on}\quad \{\chi(x)>0\}.
\end{equation*}
\end{proof}
To start with the proof of Lemma \ref{l2} let us note that for
proving (\ref{37}) it is sufficient to show that
$$
\int
(1-\psi_1)\dots(1-\psi_{r-2})\psi_{r-1}((r-1)f_\theta^{r-1}+l_{r-1})d\mu^{r-1}$$
$$+\int
(1-\psi_1)\dots(1-\psi_{r-1})\left(rf_\theta^r+v_r\right)d\mu^r
$$
\begin{equation}\label{40}
\geq\int
(1-\psi_1)\dots(1-\psi_{r-2})\left((r-1)f_\theta^{r-1}+v_{r-1}\right)d\mu^{r-1}.
\end{equation}

By Fubini's theorem the left-hand side of  (\ref{40}) is equal to
$$
\int
(1-\psi_1)\dots(1-\psi_{r-2})\psi_{r-1}((r-1)f_\theta^{r-1}+l_{r-1})d\mu^{r-1}
$$
$$
+\int (1-\psi_1)\dots(1-\psi_{r-1})\left(\int
\left(rf_\theta^r+v_r\right)d\mu(x_r)\right)d\mu^{r-1}
$$
$$=\int
(1-\psi_1)\dots(1-\psi_{r-2})[\psi_{r-1}((r-1)f_\theta^{r-1}+l_{r-1})
$$
\begin{equation}\label{41} +
(1-\psi_{r-1})\int\left(rf_\theta^r+v_r\right)d\mu(x_r)]d\mu^{r-1}.
\end{equation}

Because $f_\theta^r(x_1,\dots,x_r)$ is a joint density function of
$(X_1,\dots, X_r)$, we have
$$
\int
f_\theta^r(x_1,\dots,x_r)d\mu(x_r)=f_\theta^{r-1}(x_1,\dots,x_{r-1}),
$$
 so that the right-hand side of (\ref{41}) transforms to
$$
\int (1-\psi_1)\dots(1-\psi_{r-2})[(r-1)f_\theta^{r-1}
$$
\begin{equation}\label{42}
+\psi_{r-1}l_{r-1}+ (1-\psi_{r-1})(f_\theta^{r-1}+\int
v_rd\mu(x_r))]d\mu^{r-1}.
\end{equation}

Applying Lemma \ref{l1} with
$$
\chi=(1-\psi_1)\dots(1-\psi_{r-2}),\quad \phi=\psi_{r-1},$$
$$F_1=l_{r-1},\quad F_2=f_\theta^{r-1}+\int v_r d\mu_r,$$
we see that (\ref{42}) is greater or equal than
\begin{eqnarray}
&&
\int(1-\psi_1)\dots(1-\psi_{r-2})[(r-1)f_\theta^{r-1}+\min\{l_{r-1},f_\theta^{r-1}+\int
v_rd\mu(x_r)\}]d\mu^{r-1}\nonumber
\end{eqnarray}
\begin{equation}\label{43}
=\int (1-\psi_1)\dots(1-\psi_{r-2})[(r-1)f_\theta^{r-1}
+v_{r-1}]d\mu^{r-1},
\end{equation}
 by the definition of $v_{r-1}$ in (\ref{38}).

Moreover, by the same Lemma \ref{l1}, (\ref{42}) is equal to
(\ref{43}) if and only if (\ref{39}) is satisfied $\mu^{r-1}$-almost
anywhere on $C_{r-1}^\psi$.\end{proof}

Let now $\psi\in \Delta^N$ be any truncated stopping rule.

By (\ref{14}) we have
$$
L_N(\psi)=\sum_{n=1}^{N-1}\int
(1-\psi_1)\dots(1-\psi_{n-1})\psi_n(nf_\theta^n+l_n)d\mu^n
$$
\begin{equation}\label{43a}
+\int
(1-\psi_1)\dots(1-\psi_{N-1})\left(Nf_\theta^N+l_N\right)d\mu^r.
\end{equation}

Let $V_N^N\equiv l_N$. Applying Lemma \ref{l2} with $r=N$ and
$v_N=V_N^N$ we have
$$
L_N(\psi)\geq\sum_{n=1}^{N-2}\int
(1-\psi_1)\dots(1-\psi_{n-1})\psi_n(nf_\theta^n+l_n)d\mu^n
$$
\begin{equation}\label{44}
+\int
(1-\psi_1)\dots(1-\psi_{N-2})\left((N-1)f_\theta^{N-1}+V_{N-1}^N\right)d\mu^{N-1},
\end{equation}
where $V_{N-1}^N=\min\{l_{N-1},f_\theta^{N-1}+\int V_N^Nd\mu(x_N)
\}$.
 Also by Lemma 2, the inequality in
(\ref{44}) is in fact an equality if
\begin{equation}\label{44a}
\psi_{N-1}=I_{\{l_{N-1}\leq f_\theta^{N-1}+\int V_N^Nd\mu(x_N) \}}.
\end{equation}

 Applying  Lemma \ref{l2}  to the right-hand side of  (\ref{44}) again we see
 that
$$
L_N(\psi)\geq\sum_{n=1}^{N-3}\int
(1-\psi_1)\dots(1-\psi_{n-1})\psi_n(nf_\theta^n+l_n)d\mu^n
$$
\begin{equation}\label{45}
+\int
(1-\psi_1)\dots(1-\psi_{N-3})\left((N-2)f_\theta^{N-2}+V_{N-2}^N\right)d\mu^{N-2},
\end{equation}
where $V_{N-2}^N=\min\{l_{N-2},f_\theta^{N-2}+\int V_{N-1}^N
d\mu(x_{N-1}) \}$.  There is an equality in (\ref{45}) if
(\ref{44a}) holds and
\begin{equation}\label{46}
\psi_{N-2}=I_{\{l_{N-2}\leq f_\theta^{N-2}+\int
V_{N-1}^Nd\mu(x_{N-1}) \}},
\end{equation}
etc.

Repeating the applications of Lemma 2,
we finally get
\begin{equation}\label{47} L_N(\psi)\geq
\int \left(f_\theta^{1}+V_{1}^N\right)d\mu^1=1+\int V_1^Nd\mu(x_1),
\end{equation}
and a series of conditions on $\psi$, starting from (\ref{44a}),
(\ref{46}), etc., under which $L(\psi)$ is equal to the right-hand side of
(\ref{47}). Because Lemma 2 also
gives necessary and sufficient conditions for attaining the
equality, we also have necessary conditions for attaining the lower
bound in (\ref{47}).

In this way, formally, we have the following

\begin{theorem}\label{t3} Let $\psi\in\Delta^N$ be any (truncated) stopping rule. Then for any
$1\leq r\leq N-1$ the following inequalities hold true
$$
L_N(\psi)\geq\sum_{n=1}^{r}\int
(1-\psi_1)\dots(1-\psi_{n-1})\psi_n(nf_\theta^n+l_n)d\mu^n
$$
\begin{equation}\label{46a}
+\int
(1-\psi_1)\dots(1-\psi_{r})\left((r+1)f_\theta^{r+1}+V_{r+1}^N\right)d\mu^{r+1}
\end{equation}
$$
\geq \sum_{n=1}^{r-1}\int
(1-\psi_1)\dots(1-\psi_{n-1})\psi_n(nf_\theta^n+l_n)d\mu^n
$$
\begin{equation}\label{46b}
+\int
(1-\psi_1)\dots(1-\psi_{r-1})\left(rf_\theta^r+V_{r}^N\right)d\mu^r,
\end{equation}
where $V_N\equiv l_N$, and recursively
for $m=N-1, N-2, \dots 1 $
\begin{equation}\label{48}
V_m^N=\min\{l_m,f_\theta^m+R_m^N\},
\end{equation}
with
\begin{equation}\label{49}
   R_m^N=R_m^N(x_1,\dots, x_m)= \int
V_{m+1}^N(x_1,x_2,\dots,x_{m+1})d\mu(x_{m+1}).
\end{equation}

The lower bound in (\ref{46b}) is attained if and only if for any
$m=r,\dots, N-1$
\begin{eqnarray}
I_{\{l_{m}< f_\theta^{m}+R_m \}}\leq&\psi_{m}&\leq I_{\{l_{m}\leq
f_\theta^{m}+R_m \}} \quad \mu^m\mbox{-almost anywhere on}\quad C_m^\psi.\label{53}
\end{eqnarray}

In particular,  conditions (\ref{53}) with $m=1,2,\dots,
N-1$ are necessary and sufficient for being
$\psi=(\psi_1,\dots,\psi_{N-1},1,\dots)$ an optimal truncated rule
in $\Delta^N$. The minimum value of $L(\psi)$, over
$\psi\in\Delta^N$, is equal to
\begin{equation*}
1+\int V_1^Nd\mu(x_1)=1+R_0^N.
\end{equation*}
\end{theorem}

\begin{remark}\label{r3}\rm Despite that any $\psi=(\psi_1,\dots,\psi_{N-1},1,\dots)$ satisfying  (\ref{53}) for  $m=1,\dots, N-1$
is optimal among
all truncated  tests in $\Delta^N$, it only makes practical sense if
$$l_0> 1 +R_0^N$$
where $l_0$ defined as
$$
l_0\equiv \min_{1\leq j\leq k}\sum
_{1\leq i \leq k, i\not = j}
\lambda_{ij}.
$$

 The reason is  that $l_0$ can
be considered as "the $L(\psi)$"
function for a trivial sequential test
$(\psi_0,\phi_0)$ which, without taking
any observations, makes a decision
according to any
$\phi_0=(\phi_0^1,\dots, \phi_0^k)$
such that
$$
\phi_0^j\leq I_{\{\sum_{i\not =
j}\lambda_{ij}=l_0\}},\quad 1\leq j\leq k.
$$

In this case, there are no observations
($N(\psi_0)=0$) and it is easily seen
that
$$L(\psi_0,\phi_0)=\sum_{1\leq i,j\leq k,i\not = j}\lambda_{ij}\alpha_{ij}(\psi_0,\phi_0)=
l_0.$$ Thus, the inequality
$$l_0\leq 1+R_0^N$$
means that the trivial test $(\psi_0,\phi_0)$ is not worse than the
best truncated test in $\Delta^N$.

Because of that, we consider $V_0^N$
defined by (\ref{48}) for $m=0$, where,
by definition, $f_\theta^0=1$, as the
minimum value of $L(\psi)$ in
$\Delta^N$, in the case it is allowed
not to take any observations.
\end{remark}

\begin{remark}\label{r5}\rm It is not difficult to see from the proof of Lemma 2, that
 the problem of the optimal
 testing when the cost of the experiment is defined as
\begin{equation}\label{50a}
\int N(\psi)d\pi(\theta),
 \end{equation} with some  measure
$\pi$ (see Remark \ref{r00}), under suitable measurability
conditions, can receive essentially the
same treatment. The corresponding
optimal stopping rule in $\Delta^N$
will be defined by
\begin{equation}\label{50b}
\psi_{r}=I_{\{l_{r}\leq \int f_\theta^{r}d\pi(\theta)+\int
V_{r+1}^Nd\mu_{r+1} \}}
\end{equation}
for $r=1,2,\dots, N-1$, with $V_r^N$
defined recursively as
\begin{equation}\label{50c}
V_{r-1}^N=\min\{l_{r-1},\int f_\theta^{r-1}d\pi(\theta)+\int
V_{r}^Nd\mu(x_{r})\},
\end{equation}
starting from $r=N$, in which case
$V_N^N\equiv l_N$.

In the Bayesian context of Remark \ref{r00} the optimality of
(\ref{50b}) -- (\ref{50c}) with $\lambda_{ij}=\pi_i L_{ij}$, where
$L_{ij}$ are some non-negative losses, $i\not= j$, can be derived
also  from Theorem 5.2.2 \cite{Ghosh}. Our Theorem \ref{t3} gives,
additionally to that, a necessary condition of optimality, providing
the structure of {\em all} Bayesian truncated tests. Essentially,
they are randomizations of (\ref{50b}):
$$
I_{\{l_{r}< \int
f_\theta^{r}d\pi(\theta)+\int
V_{r+1}^Nd\mu_{r+1} \}}\leq
\psi_{r}\leq I_{\{l_{r}\leq \int
f_\theta^{r}d\pi(\theta)+\int
V_{r+1}^Nd\mu_{r+1} \}},
$$
for $r=1,2,\dots, N-1$.

In purely
Bayesian context, such conditions may
be irrelevant, because any Bayesian
test gives the same (minimum) value of
the Bayesian risk.
Nevertheless, for our (conditional)
Problems I and II, it may be important
to have a broader class of optimal
tests, for easier compliance with
(\ref{6})  in Theorem \ref{t1} (or with
(\ref{6BisBis}) in Theorem \ref{t1a}),
just like the randomization of decision
rule is important for finding tests
with a given $\alpha$-level in the
Neyman-Pearson problem (see, for
example, \cite{Lehmann}).

\end{remark}

\section{\hs GENERAL STOPPING RULES}\label{s5}

In this section we characterize the
structure of general stopping rules
minimizing $L(\psi)$.

 Let us define for any stopping rule $\psi$
$$
L_N(\psi)=\sum_{n=1}^{N-1}\int
(1-\psi_1)\dots(1-\psi_{n-1})\psi_n(nf_\theta^n+l_n)d\mu^n
$$
\begin{equation}\label{50}
+\int(1-\psi_1)\dots(1-\psi_{N-1})\left(Nf_\theta^N+l_N\right)d\mu^{N}.
\end{equation}
(cf. (\ref{43a})). This is the Lagrange-multiplier function for
$\psi$ truncated  at $N$, i.e. the rule with the components
$\psi^N=(\psi_1,\psi_2,\dots,\psi_{N-1},1,\dots)$,
$L_N(\psi)=L(\psi^N)$.

Because $\psi^N$ is truncated, the
results of the preceding section apply,
in particular, the inequalities of
Theorem \ref{t3}.

The idea of what follows is to make
$N\to\infty$, to obtain some lower
bounds for $L(\psi)$ from (\ref{46a}) -
(\ref{46b}).

To be able to do this, we need some
"approximation properties" for
$L(\psi)$, to guarantee that
$L_N(\psi)\to L(\psi)$, as
$N\to\infty$, at least for stopping
rules $\psi$ for which
$P_\theta(\tau_\psi<\infty)=1$.


\begin{lemma}\label{l3} Suppose that $\psi$ is a stopping rule such that $P_\theta(\tau_\psi<\infty)=1$.
\begin{itemize}
\item[(i)] If $L(\psi)<\infty$ and
\begin{equation}\label{59}
    \int (1-\psi_1)\dots (1-\psi_{n-1})l_nd\mu^n\to
    0,\quad\mbox{as}\quad n\to\infty,
\end{equation}
then
$$
\lim_{N\to\infty}L_N(\psi)=L(\psi).
$$
\item[(ii)] If $L(\psi)=\infty$ then $L_N(\psi)\to\infty$.
\end{itemize}
\end{lemma}

\begin{proof}
Let $L(\psi)<\infty$. Let us calculate the difference between
$L(\psi)$ and $L_N(\psi)$ in order to show that it goes to zero as
$N\to\infty$. By (\ref{50})
$$
L(\psi)-L_N(\psi)=\sum_{n=1}^\infty\int
(1-\psi_1)\dots(1-\psi_{n-1})\psi_n(nf_\theta^n+l_n)d\mu^n
$$
$$
-\sum_{n=1}^{N-1}\int
(1-\psi_1)\dots(1-\psi_{n-1})\psi_n(nf_\theta^n+l_n)d\mu^n $$
$$-\int
(1-\psi_1)\dots(1-\psi_{N-1})\left(Nf_\theta^N+l_N\right)d\mu^N
$$
$$
=\sum_{n=N}^\infty \int
(1-\psi_1)\dots(1-\psi_{n-1})\psi_n(nf_\theta^n+l_n)d\mu^n
$$
\begin{equation}\label{51}
-\int
(1-\psi_1)\dots(1-\psi_{N-1})\left(Nf_\theta^N+l_N\right)d\mu^N.
\end{equation}

The first summand converges to zero, as $N\to\infty$, being the tail
of a convergent series (this is because $L(\psi)<\infty$).

We have further
$$\int (1-\psi_1)\dots(1-\psi_{N-1})l_Nd\mu^N\to 0
$$
as $N\to\infty$, because of (\ref{59}).

It remains to show that
\begin{equation}\label{52}
\int
(1-\psi_1)\dots(1-\psi_{N-1})Nf_\theta^Nd\mu^N=NP_{\theta}(\tau_\psi\geq
N)\to 0 \;\mbox{as}\; N\to\infty.
\end{equation}

But this is again due to the fact that  $L(\psi)<\infty$ which
implies that $$E_{\theta}\tau_\psi=\sum_{n=1}^\infty
nP_{\theta}(\tau_\psi=n)<\infty.$$ Because this series is
convergent, $\sum_{n=N}^\infty nP_{\theta}(\tau_\psi=n)\to 0$. Thus,
using the Chebyshev inequality we have
$$
NP_{\theta}(\tau_\psi\geq N)\leq E_{\theta}\tau_\psi
I_{\{\tau_\psi\geq N\}}=\sum_{n=N}^\infty
nP_{\theta}(\tau_\psi=n)\to 0
$$
as $N\to\infty$, which completes the proof of (\ref{52}).

Let now $L(\psi)=\infty$.

This means that
$$\sum_{n=1}^{\infty}\int
(1-\psi_1)\dots (1-\psi_{n-1})\psi_n(nf_\theta^n+l_n)d\mu^n=\infty$$
which immediately implies by (\ref{50}) that
$$L_N(\psi)\geq\sum_{n=1}^{N-1}\int
(1-\psi_1)\dots
(1-\psi_{n-1})\psi_n(nf_\theta^n+l_n)d\mu^n\to\infty.$$
\end{proof}

Lemma \ref{l3} gives place to the following definition.

Let us say that our testing problem is {\em truncatable} if
(\ref{59}) holds for any $\psi$ with $E_\theta\tau_\psi<\infty$.

From Lemma \ref{l3} it immediately follows
\begin{corollary}\label{c1}
For any truncatable problem
$$
L_N(\psi)\to L(\psi),\quad\mbox{as}\quad N\to\infty,
$$
for any stopping rule $\psi$ such that $P_\theta(\tau_\psi<\infty)=1$.
\end{corollary}

\begin{remark}\label{r1}\rm
It is obvious from (\ref{59}) that a testing problem is
truncatable, in particular,  if
\begin{equation}\label{21}
\int l_n d\mu^n \to 0,\quad\mbox{as}\quad n\to\infty.
\end{equation}
Let us denote by $\alpha_{ij}(n,\phi)$ the error probability of a
test corresponding to a fixed number $n$ of observations, when the
decision rule $\phi$ is applied. From Theorem \ref{t2} it follows
that the left-hand side of (\ref{21}) is the minimum weighted error
sum:
$$
\int l_n d\mu^n=\inf_{\phi}\sum_{1\leq i,j\leq
k,i\not=j}\lambda_{ij}\alpha_{ij}(n,\phi),
$$
where the infimum is taken over all decision rules $\phi$.

Thus, (\ref{21}) requires a very natural behaviour of a statistical
testing problem, namely that the minimum weighted error sum, over
all fixed-sample size tests,  tend to zero, as the sample size $n$
tends to infinity.
\end{remark}
\begin{remark}\label{r6} \rm
Any Bayesian problem (with $N(\psi)=\sum_{i=1}^kN(\theta_i;\psi)\pi_i$ in (\ref{4}), where $\pi_i>0$, $i=1,\dots, k$) is truncatable.
Indeed, if $N(\psi)<\infty$ then $E_{\theta_i}\tau_\psi=N(\theta_i;\psi)<\infty$ for all $i=1,\dots,k$. Because of this,
$$
 \int (1-\psi_1)\dots (1-\psi_{n-1})l_nd\mu^n\leq  \int (1-\psi_1)\dots (1-\psi_{n-1})(\sum_{i=2}^k\lambda_{i1}f_{\theta_i}^n)d\mu^n
$$
$$
\leq\sum_{i=2}^k\lambda_{i1}P_{\theta_i}(\tau_\psi\geq n) \to 0,\quad\mbox{as}\quad n\to\infty,
$$
thus, (\ref{59}) is fulfilled.
\end{remark}
Our main results below will refer to truncatable testing problems.

To go on with the plan of passing to the limit, as $N\to\infty$, in
the inequalities of Theorem \ref{t3}, let us turn now to the
behaviour of $V_r^N$, as $N\to\infty$.

\begin{lemma}\label{l4} For any $r\geq 1$ and for any  $N\geq r$
\begin{equation}\label{52a}
V_r^N\geq V_r^{N+1}.
\end{equation}
\end{lemma}

\begin{proof} By induction over $r=N,N-1,\dots,1$.

Let $r=N$. Then by (\ref{48})
$$V_N^{N+1}=\min\{l_N,f_{\theta}^N+\int V_{N+1}^{N+1}d\mu(x_{N+1})\}\leq
l_N=V_N^N.
$$

 If we suppose that (\ref{52a}) is satisfied for some $r$, $N\geq
 r>1$, then
$$
V_{r-1}^N=\min\{l_{r-1},f_\theta^{r-1}+\int V_r^{N}d\mu(x_r)\}$$
$$\geq\min\{l_{r-1},f_\theta^{r-1}+\int V_r^{N+1}d\mu(x_r)\}=V_{r-1}^{N+1}.
$$
Thus, (\ref{52a}) is satisfied for $r-1$ as well, which completes
the induction.
\end{proof}

It follows from Lemma \ref{l4} that for any fixed $r\geq 1$ the
sequence $V_r^N$ is non-increasing. So, there exists
\begin{equation}\label{54}V_r= \lim_{N\to\infty}
V_r^N.
\end{equation}

Now, everything is prepared for passing to the limit, as
$N\to\infty$, in (\ref{46a}) and (\ref{46b}) with $\psi=\psi^N$. If $P_\theta(\tau_\psi<\infty)=1$, then the
left-hand side of (\ref{46a}) by Lemma \ref{l3} tends to $L(\psi)$,
whereas  passing to the limit in the other two parts under the
integral sign is justified by the Lebesgue monotone convergence
theorem, in view of  Lemma \ref{l4}. For the same reason, passing to
the limit as $N\to\infty$ is possible in (\ref{48}) (see
(\ref{49})).

In this way, for a truncatable testing problem we get the following

\begin{theorem}\label{t4} Let $\psi$ be any stopping rule. Then for
any $r\geq 1$ the following inequalities hold
$$
 L(\psi)\geq\sum_{n=1}^{r}\int
(1-\psi_1)\dots(1-\psi_{n-1})\psi_n(nf_\theta^n+l_n)d\mu^n
$$
\begin{equation}\label{55}
+\int(1-\psi_1)\dots(1-\psi_{r})\left((r+1)f_\theta^{r+1}+V_{r+1}\right)d\mu^{r+1}
\end{equation}
$$
\geq \sum_{n=1}^{r-1}\int
(1-\psi_1)\dots(1-\psi_{n-1})\psi_n(nf_\theta^n+l_n)d\mu^n
$$
\begin{equation}\label{56}
+\int
(1-\psi_1)\dots(1-\psi_{r-1})\left(rf_\theta^r+V_{r}\right)d\mu^r,
\end{equation}
where
\begin{equation}\label{57bis}
V_m=\min\{l_m,f_{\theta}^m+R_m\}
\end{equation}
with
\begin{equation}\label{22}
   R_m=R_m(x_1,\dots, x_m)=\int V_{m+1}(x_1,\dots,x_{m+1})d\mu(x_{m+1})
\end{equation}
 for any $m\geq 1$.

In particular,  the following lower bound holds true:
\begin{equation}\label{57a}
L(\psi)\geq 1+\int V_1 d\mu(x_1)=1+R_0.
\end{equation}
 \end{theorem}

In comparison to Theorem \ref{t3},  Theorem \ref{t4} is lacking a
very essential element: the structure of the test achieving the
lower bound on the right-hand side of (\ref{57a}). In case this test
exists, by virtue of (\ref{57a}) it has to be optimal.

First of all, let us show that if the optimal test exists, it
reaches the lower bound on the right-hand side of (\ref{57a}). More
exactly, we prove
\begin{lemma}\label{l5} For any truncatable testing problem
\begin{equation}\label{50g}
\inf_{\psi}L(\psi)=1+R_0.
\end{equation}
\end{lemma}
\begin{proof}
Let us denote
$$
U=\inf_{\psi}L(\psi),\quad U_N=1+R_0^N,
$$
where $R_0^N$ is defined in Theorem \ref{t3}.

 By Theorem \ref{t3}, for any
$N=1,2,\dots$
$$
U_N=\inf_{\psi\in\Delta^N}L(\psi).
$$
Obviously, $U_N\geq U$ for any $N=1,2,\dots$, so
\begin{equation}\label{50i}
\lim_{N\to\infty}U_N\geq U.
\end{equation}

Let us show first that in fact there is an equality in (\ref{50i}).

 Suppose the contrary, i.e. that $\lim_{N\to\infty}U_N= U+4\epsilon$, with some
 $\epsilon>0$.
 We immediately have from this that
 \begin{equation}\label{50j}U_N\geq U+3\epsilon\end{equation} for all sufficiently large $N$.

On the other hand, by the definition of $U$ there exists a $\psi$
such that $U\leq L(\psi)\leq U+\epsilon$.

Because, by Lemma \ref{l3}, $L_N(\psi)\to L(\psi)$, as $N\to\infty$,
we have that
\begin{equation}\label{50k}
L_N(\psi)\leq U+2\epsilon
\end{equation}
for all sufficiently large $N$ as well. Because, by definition,
$L_N(\psi)\geq U_N$, we have that
$$
U_N\leq U+2\epsilon
$$
for all sufficiently large $N$, which contradicts (\ref{50j}).

Thus, $$\lim_{N\to\infty}U_N=U.$$

Now, to get (\ref{50g}) we note  that, by the Lebesgue's monotone
convergence theorem,
$$
U=\lim_{N\to\infty}U_N=1+\lim_{N\to\infty}\int V_1^N(x)d\mu(x)=
1+\int V_1(x)d\mu(x)=1+R_0,
$$
thus,
$
U=1+R_0.
$
\end{proof}
\begin{remark}\label{r4}\rm
For the Bayesian context (see Remark \ref{r00}), Lemma \ref{l5} can
be derived from Theorem 5.2.3 \cite{Ghosh} if (\ref{21}) is supposed
(see also  Section 7.2 of \cite{Ferguson} or Section 9.4 of
\cite{Zacks}).
\end{remark}

The following theorem gives the structure of the optimal stopping rule
for a truncatable testing problem.

\begin{theorem}\label{t5}
\begin{equation}\label{60}
L(\psi)=\inf_{\psi^\prime} L(\psi^\prime),
\end{equation}
if and only if
\begin{equation}\label{60b}
I_{\{l_{m}< f_{\theta}^{m}+R_{m}\}}\leq\psi_{m}\leq I_{\{l_{m}\leq
f_{\theta}^{m}+R_{m}\}}\quad \mu^{m}\mbox{-almost anywhere on}\quad C_m^\psi
\end{equation}
for all $m=1, 2\dots$.
\end{theorem}

\begin{proof} Let $\psi$ be any stopping rule. By Theorem \ref{t4} for
any fixed $r\geq 1$ the following inequalities hold:
$$
L(\psi)\geq\sum_{n=1}^{r}\int
(1-\psi_1)\dots(1-\psi_{n-1})\psi_n(nf_\theta^n+l_n)d\mu^n
$$
\begin{equation}\label{61bis}
+\int
(1-\psi_1)\dots(1-\psi_{r})\left((r+1)f_\theta^{r+1}+V_{r+1}\right)d\mu({x_{r+1}})
\end{equation}
$$
\geq \sum_{n=1}^{r-1}\int
(1-\psi_1)\dots(1-\psi_{n-1})\psi_n(nf_\theta^n+l_n)d\mu^n
$$
\begin{equation}\label{62bis}
+\int
(1-\psi_1)\dots(1-\psi_{r-1})\left(rf_\theta^r+V_{r}\right)d\mu(x_{r})
\end{equation}
$$
\geq \dots
$$
\begin{equation}\label{63bis}
\geq \int\psi_1(f_\theta^1+l_1)d\mu^1+\int
(1-\psi_1)\left(2f_\theta^{2}+V_{2}\right)d\mu(x_2)
\end{equation}
\begin{equation}\label{64bis}
\geq 1+\int V_{1}d\mu(x_1)=1+R_0.
\end{equation}

Let us suppose that $L(\psi)=1+R_0$. Then, by Lemma \ref{l5}, there are equalities in
all the inequalities (\ref{61bis})-(\ref{64bis}). Applying the "only
if"-part of Lemma \ref{l2}  and using (\ref{57bis}) and (\ref{22}),
successively, starting from the last inequality (\ref{64bis}), we
get that (\ref{60b}) has to be satisfied for any $m=1,2,\dots$. The
first part of the Theorem is proved.

Let now $\psi$ be any test satisfying (\ref{60b}).

Applying the "if"-part of Lemma \ref{l2}  and using (\ref{57bis})
and (\ref{22}) again, we see that all the inequalities in
(\ref{62bis})-(\ref{64bis}) are in fact equalities for
$$\psi^{r}=(\psi_1,\psi_2,\dots,\psi_{r},1,\dots).$$

In particular, this means that there exists
$$
\lim_{r\to\infty}[\sum_{n=1}^{r}\int
(1-\psi_1)\dots(1-\psi_{n-1})\psi_n(nf_\theta^n+l_n)d\mu^n
$$
\begin{equation}\label{65bis}
+\int
(1-\psi_1)\dots(1-\psi_{r})\left((r+1)f_\theta^{r+1}+V_{r+1}\right)d\mu(x_{r+1})]=1
+R_0.
\end{equation}

It follows from (\ref{65bis}) that
\begin{eqnarray*}\
&\limsup_{r\to\infty}\int
(1-\psi_1)\dots(1-\psi_{r})(r+1)f_\theta^{r+1}d\mu(x_{r+1})\\
&=\limsup_{r\to\infty}(r+1)P_\theta(\tau_\psi\geq r+1)
\leq 1+R_0,
\end{eqnarray*}
which implies that $\lim_{r\to\infty}P_\theta(\tau_\psi\geq r+1)=0$. Thus, $P_\theta(\tau_\psi<\infty)=1$.

From (\ref{65bis}), it  follows as well that
\begin{equation}\label{66bis}
\lim_{r\to\infty}\sum_{n=1}^{r}\int[
(1-\psi_1)\dots(1-\psi_{n-1})\psi_n(nf_\theta^n+l_n)d\mu^n\leq
1+R_0.
\end{equation}

But the left-hand side of (\ref{66bis}) is $L(\psi)$ (because
$P_\theta(\tau_\psi<\infty)=1$) and hence
\begin{equation}\label{67bis}
L(\psi)\leq 1+R_0.
\end{equation}

On the other hand, by virtue of  Theorem \ref{t4},
$$
L(\psi)\geq 1+R_0,
$$
which proves, together  with  (\ref{67bis}), that
$
L(\psi)=1+R_0.
$
\end{proof}
\begin{remark}\label{r7}\rm
Once again (see Remark \ref{r3}), the optimal stopping rule $\psi$
from Theorem \ref{t5} only makes practical sense if $l_0>1+R_0$.
\end{remark}
\begin{remark}\label{r7a}\rm
From the results of this section it is not difficult to see that the
same method works as well for minimizing
$$
\int N(\psi)d\pi(\theta)
$$
(see Remark \ref{r5}).

Repeating the steps which led us to  Theorem \ref{t5} we get that
the corresponding optimal stopping rule has the form
\begin{equation}\label{90}
\psi_{r}=I_{\left\{l_{r}\leq \displaystyle{ \int}
f_\theta^{r}d\pi(\theta)+\displaystyle{\int V_{r+1}d\mu(x_{r+1})}
\right\}},\; r=1,2,3,\dots,
\end{equation}
with
$$
V_r=\lim_{N\to\infty}V_r^N,
$$
being $V_r^N$ defined for $r=N-1,N-2,\dots, 1$ recursively by
$$
V_{r}^N=\min\{l_r, \int f_\theta^{r}d\pi(\theta)+\displaystyle{\int
V_{r+1}^Nd\mu(x_{r+1})} \}
$$
starting from $V_N^N\equiv l_N$.

In a particular case of Remark \ref{r00} and
$$
\int N(\theta;\psi)d\pi(\theta)=\sum_{i=1}^k\pi_i N(\theta_i;\psi)
$$
being $\lambda_{ij}=L_{ij}\pi_i$,  this gives an optimal stopping
rule for the Bayesian problem considered in \cite{Cochlar}.

In particular, for $k=2$, this gives an optimal stopping rule for
the Bayesian problem considered in \cite{CochlarVrana}.
\end{remark}

\normalsize
\section{\hs APPLICATIONS TO THE CONDITIONAL PROBLEMS}\label{s7}

In this section, we apply the results obtained in the preceding
sections to minimizing the average sample size $N(\psi)=E_{\theta}
\tau_\psi$ over all sequential testing procedures  with error
probabilities not exceeding some prescribed levels (see Problems I
and II in Section 1). Recall that we are supposing that our problems
are truncatable (see Section 5).

Combining Theorems \ref{t1}, \ref{t2} and  \ref{t5}, we immediately
have the following solution to Problem I.
\begin{theorem}\label{t8} Let $\psi$ satisfy
(\ref{60b}) for all $m=1,2,\dots$, with any $\lambda_{ij}>0$, $i,j=1,\dots,k$, $i\not =j$,
(recall that $R_m$ and $l_m$ in (\ref{60b}) are functions of $\lambda_{ij}$), and let $\phi$ be any
decision rule satisfying (\ref{7aa}).

Then for all sequential testing procedures
$(\psi^\prime,\phi^\prime)$ such that
\begin{equation}\label{7.2}
    \alpha_{ij}(\psi^\prime,\phi^\prime)\leq
    \alpha_{ij}(\psi,\phi)\quad \mbox{for all}\quad
    i,j=1,\dots, k,\; i\not=j,
\end{equation}
it holds
\begin{equation}\label{7.3}
    N(\psi^\prime)\geq N(\psi).
\end{equation}

The inequality in (\ref{7.3}) is strict if at least one of the
inequalities in (\ref{7.2}) is strict.

If there are equalities in all of the  inequalities in (\ref{7.2})
and (\ref{7.3}), then $\psi^\prime$ satisfies (\ref{60b}) for all $m=1,2,\dots$ as well
(with $\psi^\prime$ instead of   $\psi$).
\end{theorem}
\begin{proof} The only thing  to be proved is the last
assertion.

Let us suppose that
$$\alpha_{ij}(\psi^\prime,\phi^\prime)=
    \alpha_{ij}(\psi,\phi), \quad \mbox{for all}\quad
    i,j=1,\dots, k,\; i\not=j,$$
     and $$N(\psi^\prime)=
    N(\psi).$$
Then, obviously,
\begin{equation}\label{7.4}
    L(\psi^\prime,\phi^\prime)=L(\psi,\phi)=L(\psi)\geq
    L(\psi^\prime)
\end{equation}
(see (\ref{4}) and Remark \ref{r2}).

By Theorem \ref{t5}, there can not be strict inequality  in (\ref{7.4}), so $L(\psi)=L(\psi^\prime)$.
From Theorem \ref{t5} it follows now that $\psi^\prime$ satisfies
(\ref{60b}) as well.
\end{proof}

Analogously, combining Theorems \ref{t1a}, \ref{t2} and  \ref{t5},
we also have the following solution to Problem II.
\begin{theorem}\label{t9} Let $\psi$ satisfy
(\ref{60b}) for all $m=1,2,\dots$, with $\lambda_{ij}=\lambda_i$ for all $j=1,\dots, k$, where $\lambda_i>0$, $i=1,\dots k$ are any numbers,
 and let $\phi$ be any decision rule
such that $$ \phi_{nj}\leq I_{\left\{\sum_{i\not
=j}\lambda_{i}f_{\theta_i}^n=\min_j\sum_{i\not
=j}\lambda_{i}f_{\theta_i}^n\right\}}$$ for all $j=1,\dots, k$ and
for all $n=1,2,\dots$.

Then for any sequential test $(\psi^\prime,\phi^\prime)$ such that
\begin{equation}\label{7.102}
    \beta_{i}(\psi^\prime,\phi^\prime)\leq
    \beta_{i}(\psi,\phi)\quad \mbox{for any}\quad
    i=1,\dots, k,
\end{equation}
it holds
\begin{equation}\label{7.103}
    N(\psi^\prime)\geq N(\psi).
\end{equation}

The inequality in (\ref{7.103}) is strict if at least one of the
inequalities in (\ref{7.102}) is strict.

If there are equalities in all of the  inequalities in (\ref{7.102})
and (\ref{7.103}), then $\psi^\prime$ satisfies
(\ref{60b}) for all $m=1,2,\dots$ as well (with $\psi^\prime$ instead of   $\psi$).
\end{theorem}
\begin{remark}\label{r9}\rm
There are examples of applications of Theorem \ref{t8} (or
\ref{t9}), in the case of two simple hypotheses  based on
independent observations, in \cite{NovikovIJPAM}.

A numerical example related to the modified Kiefer-Weiss problem for
independent and identically distributed observations can be found in
\cite{Lorden}. Obviously, our Theorem \ref{t8} provides, for this
particular case, randomized versions of the optimal
sequential test studied in \cite{Lorden} (see also
\cite{NovikovIJPAM}).
\end{remark}
 \vspace{4mm}
\section*{ACKNOWLEDGEMENTS}
\small
 The author greatly appreciates the support of the Autonomous
 Metropolitan University, Mexico City, Mexico, where this work was
 done, and the support of the National System of Investigators (SNI)
 of
 CONACyT, Mexico.

This work is also partially supported by Mexico's CONACyT Grant no.
CB-2005-C01-49854-F.

The author thanks the anonymous referees for reading the article
carefully and for their valuable suggestions and comments. \vspace{3mm}
 \footnotesize
\begin{flushright}
(Received November 2, 2007.)\,\ \rule{0mm}{0mm}
\end{flushright}

\vspace*{2mm}

{\mi
\begin{flushright}
\begin{minipage}[]{124mm}
{Andrey Novikov, Departamento de Matem\'aticas, Universidad
Aut\'onoma Metropolitana - Unidad Iztapalapa, San Rafael Atlixco
186, col. Vicentina, C.P. 09340, M\'exico D.F., M\'exico
\\ e-mail: {\tt an@xanum.uam.mx}
\\{\tt http://mat.izt.uam.mx/profs/anovikov/en}}
\end{minipage}
\end{flushright}
}

\end{document}